\documentclass[10pt]{article}
\font\smallit=cmti10
\font\smalltt=cmtt10

\usepackage{amssymb,latexsym,amsmath,epsfig,amsthm} 

\makeatletter

\renewcommand\section{\@startsection {section}{1}{\z@}
{-30pt \@plus -1ex \@minus -.2ex}
{2.3ex \@plus.2ex}
{\normalfont\normalsize\bfseries}}

\renewcommand\subsection{\@startsection{subsection}{2}{\z@}
{-3.25ex\@plus -1ex \@minus -.2ex}
{1.5ex \@plus .2ex}
{\normalfont\normalsize\bfseries}}

\renewcommand{\@seccntformat}[1]{\csname the#1\endcsname. }

\makeatother

\newtheorem{theorem}{Theorem}

\newtheorem{corollary}[theorem]{Corollary}
\newtheorem{proposition}[theorem]{Proposition}

\newtheorem{remark}[theorem]{Remark}

\begin{document}

\begin{center}
\uppercase{\bf On the Stern sequence and its twisted version}
\vskip 20pt
{\bf J.-P. Allouche}\\
{\smallit CNRS, Institut de Math., \'Equipe Combinatoire et Optimisation \\
Universit\'e Pierre et Marie Curie, Case 247 \\
4 Place Jussieu \\
F-75252 Paris Cedex 05 \\
France \\
}
{\tt allouche@math.jussieu.fr}\\
\end{center}
\vskip 30pt

\centerline{\smallit Received: , Revised: , Accepted: , Published: }
\vskip 30pt

\centerline{\bf Abstract}

\noindent
In a recent preprint on ArXiv, Bacher introduced a twisted version of the Stern 
sequence. His paper contains in particular three conjectures relating the generating series
for the Stern sequence and for the twisted Stern sequence. Soon afterwards Coons published
two papers in {\it Integers}: first he proved these conjectures, second he used his result
to obtain a correlation-type identity for the Stern sequence. We recall here a simple
result of Reznick and we state a similar result for the twisted Stern sequence.
We deduce an easy proof of Coons' identity, and a simple proof of Bacher's conjectures.
Furthermore we prove identities similar to Coons' for variations on the Stern sequence 
that include Bacher's sequence.

\pagestyle{myheadings}
\markright{\smalltt INTEGERS: xx (XXXX)\hfill}
\thispagestyle{empty}
\baselineskip=12.875pt
\setcounter{page}{1} 
\vskip 30pt

\section{Introduction}

The Stern sequence is a sequence of integers ${\mathbf s} = (s(n))_{n \geq 0}$
that can be defined inductively by $s(0) = 0$, $s(1) = 1$, and for all $n \geq 1$,
$s(2n) = s(n)$ and $s(2n+1) = s(n)+s(n+1)$. (Note that these two equalities are 
actually true for all $n \geq 0$.) This is sequence~A002487 in \cite{oeis}. 
Its first few terms are 
$$
0, 1, 1, 2, 1, 3, 2, 3, 1, 4, 3, 5, 2, 5, 3, 4, 1 \ldots
$$
Several authors studied that sequence, see, e.g., \cite{Urbiha, Northshield} and the 
references therein. (Note that some authors call Stern sequence the shifted sequence 
$(s(n+1))_{n \geq 0}$.)

\medskip

Bacher introduced recently in \cite{Bacher} a twisted version of the Stern sequence
${\mathbf t} = (t_n)_{n \geq 0}$ defined inductively by $t(0) = 0$, $t(1) = 1$, and 
for all $n \geq 1$, $t(2n) = -t(n)$, $t(2n+1)= -t(n) - t(n+1)$. He gave several 
interesting properties of the sequences ${\mathbf s}$ and ${\mathbf t}$ and formulated
conjectural relations between the generating series $\sum s(n) X^n$, $\sum t(3.2^e + n) X^n$, 
$\sum (s(2+n)-s(1+n))X^n$, and $\sum(t(2+n)+t(1+n))X^n$.

\medskip

In the recent paper \cite{Coons-conj} Coons proved Bacher's conjectures. He then used
in \cite{Coons} his results to prove the following identity for the Stern sequence: 
if $e$ and $r$ are integers with $e \geq 0$, then for every integer $n \geq 0$
$$
s(r)s(2n+5) + s(2^e-r)s(2n+3) = s(2^e(n+2)+r) + s(2^e(n+1)+r).
$$

\medskip

We recall here (see Section~\ref{Rez}) a result of Reznick in \cite{Reznick}, and we
deduce an easy proof of Coons' identity. We also prove a result similar to Reznick's result
for the Bacher-Stern sequence which yields a short proof of Bacher's conjectures.
Furthermore we prove identities analogous to Reznick's and Coons' identities for sequences 
satisfying recurrence relations similar to Stern's which include Bacher's sequence.

\section{Three auxiliary results}

We start with three propositions. The first one is \cite[Corollary~4]{Reznick} for 
which Reznick gives a short proof.

\begin{proposition}{\rm \cite{Reznick}}\label{main}
Let $e$ and $r$ be integers with $e \geq 0$ and $0 \leq r \leq 2^e$.
Then, for every integer $n \geq 0$, we have
$$
s(2^e n + r) = s(r)s(n+1) + s(2^e - r)s(n).
$$
\end{proposition}

The next Proposition is similar to Proposition~\ref{main}.

\begin{proposition}\label{sim-main}
Let $e$ and $r$ be integers with $e \geq 0$ and $0 \leq r \leq 2^e$.
Then, for every integer $n \geq 1$, we have
$$
t(2^e n + r) = (-1)^e(s(r)t(n+1) + s(2^e - r)t(n)).
$$
\end{proposition}

\proof We prove by induction on $e \geq 0$ that, for every $r \in [0, 2^e]$, the identity
in the proposition holds. This is immediate for $e=0$ (thus $r \in \{0, 1\}$). If the result
is true for some $e$, then, using the definition of ${\mathbf t}$, the induction hypothesis,
and the definition of ${\mathbf s}$, we have

\begin{itemize}

\item if $2r \in [0, 2^{e+1}]$, then
$$
\begin{array}{ll}
t(2^{e+1} n + 2r) &= - t(2^e n + r) = (-1)^{e+1} (s(r)t(n+1) + s(2^e - r)t(n)) \\
&= (-1)^{e+1}(s(2r)t(n+1) + s(2^{e+1} - 2r)t(n)); \\
\end{array}
$$

\item if $2r+1 \in [0, 2^{e+1}]$, then
$$
\begin{array}{ll}
t(2^{e+1} n + 2r + 1) &= t(2(2^e n + r)+1) = - t(2^e n + r) - t(2^e n + r + 1) \\
&= \left\{ 
\begin{array}{ll}
        (-1)^{e+1}(s(r)t(n+1) + s(2^e - r)t(n)) \\
         \ \ + (-1)^{e+1}(s(r+1)t(n+1) + s(2^e - r - 1)t(n)) \\
\end{array}
\right. \\
&= \left\{ 
\begin{array}{ll}
        (-1)^{e+1} (s(r) + s(r+1))t(n+1) \\
         \ \ + (-1)^{e+1} (s(2^e - r) + s(2^e - r - 1))t(n) \\
\end{array}
\right. \\
&= (-1)^{e+1} (s(2r+1)t(n+1) + s(2(2^e - r - 1) + 1) \\
&= (-1)^{e+1} (s(2r+1)t(n+1) + s(2^{e+1} - 2r - 1)t(n).
\end{array}
$$

\end{itemize}

The last result we need is a consequence of Proposition~\ref{main}.

\begin{proposition}\label{sum-s}
Let $S(X) = \sum_{n \geq 0} s(N) X^n$. Then
$$
S(X) = S(X^{2^e}) \sum_{0 \leq r \leq 2^e-1} (s(2^e-r)X^r + s(r)X^{r-2^e}).
$$
\end{proposition}

\proof This is an easy consequence of Proposition~\ref{main} 
(also recall that $s(0)=0$): we write
$$
\begin{array}{ll}
S(X) &= \displaystyle\sum_{n \geq 0} s(n) X^n = 
\displaystyle\sum_{0 \leq r \leq 2^e-1} \sum_{k \geq 0} s(k.2^e + r) X^{k.2^e + r} \\
&= \displaystyle\sum_{0 \leq r \leq 2^e-1} X^r 
     \sum_{k \geq 0}(s(r)s(k+1) + s(2^e - r)s(k)) X^{k.2^e} \\
&= \displaystyle\sum_{0 \leq r \leq 2^e-1} (s(r) X^{r-2^e} + s(2^e - r)X^r)
     \sum_{k \geq 0} s(k)X^{k.2^e} \\
&= S(X^{2^e}) \displaystyle\sum_{0 \leq r \leq 2^e-1} (s(r) X^{r-2^e} + s(2^e - r)X^r).
\end{array}
$$

\begin{remark}
{\rm Note that, as indicated by the referee, Proposition~\ref{sum-s} is essentially Lemma~8 of 
\cite{Coons-conj} which states that for all $k \geq 0$
$$
X \prod_{0 \leq i \leq k-1} \left(1 + X^{2^i} + X^{2^{i+1}}\right)
= \sum_{1 \leq n \leq 2^k} s(n) X^n + \sum_{1 \leq n \leq 2^k-1} s(2^k-n) X^{n+2^k}.
$$
This can be deduced from the following property of the generating function for $(s(n))$
$$
S(X^2) = \left(\frac{X}{1+X+X^2}\right) S(X)
$$ 
(see the paper of Carlitz \cite[p.\ 19]{Carlitz}, where the shifted sequence
$(\Theta_0(n)) = (s(n+1))$ is studied).}
\end{remark}

\section{A direct proof of Coons' identity}\label{Rez}

Theorem~1 of \cite{Coons} is a straightforward corollary of Reznick's result 
(Proposition~\ref{main} above).

\begin{corollary}\label{Coons}
Let $e$ and $r$ be integers with $e \geq 0$ and $0 \leq r \leq 2^e$.
Then, for every integer $n \geq 0$, we have
$$
s(r)s(2n+5) + s(2^e-r)s(2n+3) = s(2^e(n+2)+r) + s(2^e(n+1)+r).
$$
\end{corollary}

\proof Let $S(e,r,n) = s(2^e(n+2)+r) + s(2^e(n+1)+r)$. Applying Proposition~\ref{main} 
with $n$ replaced by $n+2$ and $n+1$, and the definition of the sequence 
${\mathbf s}$ yields
$$
\begin{array}{lll}
S(e,r,n)
&=& s(r)s(n+3) + s(2^e - r)s(n+2) + s(r)s(n+2) + s(2^e - r)s(n+1) \\
&=& s(r)(s(n+3)+s(n+2)) + s(2^e - r)(s(n+2)+s(n+1))\\
&=& s(r)s(2n+5) +  s(2^e - r)s(2n+3). \\
\end{array}
$$

\section{A simple proof of Bacher's conjectures}

We can now prove the three conjectures that Bacher proposed in \cite{Bacher} (Conjectures~1.3,
3.2~(i), and 3.2~(ii)) as Theorems~\ref{B-conj1}, \ref{B-conj2}, and \ref{B-conj3} below.

\begin{theorem}\label{B-conj1}
Let $S(X) = \sum_{n \geq 0} s(n) X^n$ and $T(X) = \sum_{n \geq 0} t(n) X^n$ be the 
generating series of ${\mathbf s}$ and ${\mathbf t}$. Then, there exists a series 
$U(X) = \sum_{n \geq 0} u(n) X^n$ with integral coefficients, such that
$$
\forall e \geq 0, \ \sum_{n \geq 0} t(3.2^e + n) X^n = (-1)^e U(X^{2^e}) S(X).
$$
\end{theorem}

\proof The series $U(X)$ must satisfy in particular
$\sum_{n \geq 0} t(3+n) X^n = U(X) S(X)$. This relation {\it defines} a series
$U(X)$ that clearly has integer coefficients ($s(1) = 1$, and $s(0)=0$).
Now, using Proposition~\ref{sim-main} above, the definition of $U(X)$, and 
Proposition~\ref{sum-s}, we have
$$
\begin{array}{ll}
\displaystyle\sum_{n \geq 0} t(3.2^e + n) X^n &= \
\displaystyle\sum_{0 \leq r \leq 2^e-1} \sum_{k \geq 0} t(3.2^e + k.2^e + r) X^{k.2^e + r} \\
&= \displaystyle\sum_{0 \leq r \leq 2^e-1} X^r \sum_{k \geq 0} t(2^e(3+k) + r) X^{k.2^e} \\
&= \displaystyle\sum_{0 \leq r \leq 2^e-1} X^r \sum_{k \geq 0}
    (-1)^e(s(r)t(k+4) + s(2^e - r)t(k+3)) X^{k.2^e} \\
&= \displaystyle(-1)^e\sum_{0 \leq r \leq 2^e-1} (s(2^e-r)X^r + s(r)X^{r-2^e})
   \sum_{k \geq 0} t(3+k) X^{k.2^e} \\
&= \displaystyle(-1)^e\sum_{0 \leq r \leq 2^e-1} (s(2^e-r)X^r + s(r)X^{r-2^e})
   U(X^{2^e})S(X^{2^e}) \\
&= (-1)^e S(X) U(X^{2^e}). \\
\end{array}
$$

\begin{theorem}\label{B-conj2}
Let $A(X) = \displaystyle\frac{1}{S(X)} \sum_{n \geq 0} (s(2+n) - s(1+n)) X^n$.
Then
$$
\sum_{n \geq 0} (s(2^{e+1} + n) - s(2^e + n)) X^n = A(X^{2^e}) S(X).
$$
\end{theorem}

\proof Let $A_e(X) = \displaystyle\sum_{n \geq 0} (s(2^{e+1} + n) - s(2^e + n)) X^n$.
We write, using Proposition~\ref{main} and Proposition~\ref{sum-s} (recall that $s(2)-s(1)=0$),
$$
\begin{array}{ll}
A_e(X) &=
\displaystyle\sum_{0 \leq r \leq 2^e - 1} 
\sum_{k \geq 0} (s(2^{e+1} + k.2^e + r) - s(2^e + k.2^e + r)) X^{k.2^e} \\
&= \displaystyle\sum_{0 \leq r \leq 2^e - 1} X^r
\sum_{k \geq 0} (s(2^e(k+2)+r) - s(2^e(k+1)+r)) X^{k.2^e}. \\
\end{array}
$$
Thus
$$
\begin{array}{ll}
A_e(X) &= \left\{
\begin{array}{ll}
\displaystyle\sum_{0 \leq r \leq 2^e - 1} X^r
\sum_{k \geq 0} ((s(r)s(k+3) + s(2^e-r)s(k+2)) X^{k.2^e} \\
- \displaystyle\sum_{0 \leq r \leq 2^e - 1} X^r (s(r)s(k+2) + s(2^e-r)s(k+1))X^{k.2^e} \\
\end{array}
\right. \\
&= \left\{
\begin{array}{ll}
\displaystyle\sum_{0 \leq r \leq 2^e - 1} s(r) X^r \sum_{k \geq 0} (s(k+3)-s(k+2)) X^{k.2^e} \\
+ \displaystyle\sum_{0 \leq r \leq 2^e - 1} s(2^e-r) X^r \sum_{k \geq 0} (s(k+2)-s(k+1))X^{k.2^e} \\
\end{array}
\right. \\
&= 
\displaystyle\sum_{0 \leq r \leq 2^e - 1} (s(r) X^{r-2^e} + s(2^e-r)) X^r
\displaystyle\sum_{k \geq 0} (s(k+2)-s(k+1))X^{k.2^e} \\
&= 
\displaystyle\frac{S(X)}{S(X^{2^e})} \sum_{k \geq 0} (s(k+2)-s(k+1))X^{k.2^e} = S(X) A(X^{2^e}). \\
\end{array}
$$

\begin{theorem}\label{B-conj3}
Let $B(X) = \displaystyle\frac{1}{S(X)} \sum_{n \geq 0} (t(2+n) + t(1+n)) X^n$.
Then
$$
(-1)^{e+1} \sum_{n \geq 0} (t(2^{e+1} + n) + t(2^e + n)) X^n = B(X^{2^e}) S(X).
$$
\end{theorem}

\proof The proof is the same as the proof of Theorem~\ref{B-conj2}, except that we use
Propositions~\ref{sim-main} and \ref{sum-s} instead of Propositions~\ref{main} and \ref{sum-s}.

\section{Similar sequences}

Proposition~\ref{sim-main} gives an expression of $t(2^e n + r)$ in terms of $t(n)$ and
$t(n+1)$ with coefficients in terms of ${\mathbf s}$. One might want to find relations
of the same kind but involving ${\mathbf t}$ only. In this section we give such a relation.
More generally we prove such relations for sequences satisfying recurrence relations similar
to the recurrences defining the Stern sequence.

\begin{theorem}\label{general}
Let ${\mathbf v} = (v(n))_{n \geq 0}$ be a sequence of real numbers satisfying 
$$
\exists n_0 \geq 0, \ \exists (a, b, c) \in {\mathbb R}^3, \ \forall n \geq n_0 \
v(2n) = a v(n) \ \mbox{\rm and} \ v(2n+1) = bv(n) + cv(n+1).
$$
Then, for all integers $(e,r)$ with $e \geq 0$ and $r \in [0, 2^e]$, there exist 
$A = A(e,r)$ and $B = B(e,r)$ such that for all $n \geq n_0$ 
$$
v(2^e n + r) = A(e,r) v(n) + B(e,r) v(n+1).
$$
\end{theorem}

\proof We prove by induction on $e$ that for all $r \in [0, 2^e]$, there exist 
$A(e,r)$ and $B(e,r)$ satisfying the conditions in the theorem.
For $e=0$, hence $r \in \{0,1\}$ one gets from the definition of ${\mathbf v}$ that
$A(0,0) = 1$, $B(0,0) = 0$, $A(0,1) = 0$, and $B(0,1) = 1$. Going from $e$ to $e+1$
yields $A(e+1,2r) = a A(e,r)$, $B(e+1, 2r) = a B(e,r)$, if $0 \leq 2r \leq 2^{e+1}$,
and $A(e+1,2r+1) = bA(e,r) + cA(e,r+1)$, $B(e+1, 2r+1) = bB(e,r) + cB(e,r+1)$, if
$0 \leq 2r+1 \leq 2^{e+1}$.

\begin{corollary}\label{cor-general}
Let ${\mathbf v} = (v(n))_{n \geq 0}$ be a sequence of real numbers satisfying
$$
\exists n_0 \geq 0, \ \exists (a, b, c) \in {\mathbb R}^3, \ \forall n \geq n_0 \
v(2n) = a v(n) \ \mbox{\rm and} \ v(2n+1) = bv(n) + cv(n+1).
$$
Then, for all integers $(e,r)$ with $e \geq 0$ and $r \in [0, 2^e]$, there exist
$A = A(e,r)$ and $B = B(e,r)$ such that for all $n \geq n_0$
$$
A(e,r)v(2n+3) + B(e,r)v(2n+5) = cv(2^e(n+2)+r) + bv(2^e(n+1)+r).
$$
\end{corollary}

\proof Apply Theorem~\ref{general} with $n$ replaced by $n+2$ and $n+1$
to the left side of the identity to be proven.

\begin{remark}
{\rm The quantities $A(e,r)$ and $B(e,r)$ can of course be computed in terms of $e,r$ 
and of certain values of ${\mathbf v}$. For example if the sequence ${\mathbf v}$ is not 
trivial, there exist two integers $x_0$ and $y_0$ with $x_0, y_0 \geq n_0$ such that
$
\left|
\begin{array}{cc}
v(x_0) & v(x_0+1) \\
v(y_0) & v(y_0+1) \\
\end{array}
\right| \neq 0.
$
Then 
$$
\begin{array}{ll}
v(2^e x_0 + r) &= A(e,r) v(x_0) + B(e,r) v(x_0+1) \\
v(2^e y_0 + r) &= A(e,r) v(y_0) + B(e,r) v(y_0+1) \\
\end{array}
$$
yields
$$
A(e,r) = (v(y_0)v(x_0 +1) - v(x_0)v(y_0 +1))^{-1}
(v(x_0+1)v(2^e y_0+r) - v(y_0+1)v(2^e x_0+r))
$$
and
$$
B(e,r) = (v(x_0)v(y_0 +1) - v(y_0)v(x_0 +1))^{-1}
(v(x_0)v(2^e y_0+r) - v(y_0)v(2^e x_0+r)).
$$
}
\end{remark}

\section{Examples}

\subsection{The Stern sequence again}

One can apply Theorem~\ref{general} to the Stern sequence, for which
$n_0 = 0$, $a=b=c=1$. The values of $A$ and $B$ can be obtained by taking 
$n=0$ and $n=1$ in the relation $s(2^e n + r) = A(e, r)s(n) + B(e, r)s(n + 1)$,
yielding $B(e,r) = s(r)$ and $A(e,r) = s(2^e n + r) - s(r)$. To obtain the result 
of Proposition~\ref{main} and Corollary~\ref{Coons} this way, it remains to prove that
for all $e \geq 0$ and $r \in [0, 2^e]$ one has $s(2^e + r) - s(r) = s(2^e-r)$.
This last equality can be proven by induction on $e$, but this is also
Corollary~3.1 in \cite{DS} (see also \cite[Theorem~1.2]{Bacher} where the
author adds that this identity ``is probably well-known to the experts'').

\subsection{The case of Bacher's twisted Stern sequence}

The definition of Bacher's twisted Stern sequence ${\mathbf t} = (t(n))_{n \geq 0}$
recalled in the Introduction shows that ${\mathbf t}$ satisfies the hypotheses of
Theorem~\ref{general} with $a=b=c=-1$, and $n_0 = 1$. Note that the first few terms 
of ${\mathbf t}$ are: 
$$
0, 1, -1, 0, 1, 1, 0, -1, -1, -2, -1, -1, 0, 1, 1, 2, \ldots
$$
Applying Theorem~\ref{general} and Corollary~\ref{cor-general} we get 
the following results.

\begin{theorem}\label{similar}
Let $e$ and $r$ be integers with $e \geq 0$ and $0 \leq r \leq 2^e$. Then, for every
integer $n \geq 1$, we have
$$
t(2^e n + r) = - t(2^{e+1} + r)t(n) - t(3.2^e - r)t(n+1).
$$
\end{theorem}

\proof From Theorem~\ref{general} we have the existence of $A'$ and $B'$ such that
$t(2^e n + r) = A'(e,r)t(n) + B'(e,r) t(n+1)$ for $n \geq 1$. Taking $n=2$ and using
that $t(2) = -1$ and $t(3) = 0$, we get $A'(e,r) = - t(2^{e+1}+r)$. Now taking $n=1$
yields $t(2^e + r) = A'(e,r) - B'(e,r)$. Hence $B'(e,r) = A'(e,r) - t(2^e + r)$, i.e.,
$B'(e,r) = - t(2^{e+1}+r) - t(2^e + r)$. An immediate induction on $e$ shows that for
$r \in [0, 2^e]$ one has $t(2^{e+1}+r) + t(2^e + r) = t(3.2^e - r)$. Hence the result.

\begin{corollary}
Let $e$ and $r$ be integers with $e \geq 0$ and $0 \leq r \leq 2^e$.
Then, for every integer $n \geq 0$, we have
$$
t(2^{e+1} + r)t(2n+3) + t(3.2^e-r)t(2n+5) = t(2^e(n+2)+r) + t(2^e(n+1)+r).
$$
\end{corollary}

\subsection{Other variations on Stern's sequence}

Let the three sequences $(z_1(n))_{n \geq 0}$, $(z_2(n))_{n \geq 0}$, 
and $(z_3(n))_{n \geq 0}$ defined by (using the notation of \cite{oeis}): 
for all $n \geq 0$, $z_1(n) = A005590(n)$, and for all $n \geq 1$,
$z_2(n) = A177219(n)$, and $z_3(n) = A049347(n)$ with $z_2(0) = z_3(0) = 0$.
These sequences satisfy respectively
$$
(z_1(0), z_1(1)) = (0, 1), \ \mbox{\rm and} \ \forall n \geq 1, 
z_1(2n) = z_1(n), \ z_1(2n+1) = -z_1(n) + z_1(n+1),
$$
$$
(z_2(0), z_1(1)) = (0,1), \ \mbox{\rm and} \ \forall n \geq 1, 
z_2(2n) = -z_2(n), \ z_2(2n+1) = -z_2(n) + z_2(n+1), 
$$
$$
(z_3(0), z_3(1)) = (0,1), \ \mbox{\rm and} \ \forall n \geq 1, 
z_3(2n) = -z_3(n), \ z_3(2n+1) = z_3(n) + z_3(n+1). 
$$
Note that he last sequence $(z_3(n))_{n \geq 0}$ is the $3$-periodic sequence 
with period $(0, 1, -1)$ (hint: prove by induction on $n$ that for all $j \leq n$ 
one has $(z_3(3j), z_3(3j+1),z_3(3j+2)) = (0, 1, -1)$). Also note that all relations 
$z_i(2n) = \pm z_i(n)$ and $z_i(2n+1) = \pm z_i(n) + z_i(n+1)$, $i = 1, 2, 3$, are
actually valid for $n \geq 0$.

We know from Theorem~\ref{general} that, for all $e \geq 0$ and $r \in [0, 2^e]$,
there exist $A_i(e,r)$ and $B_i(e,r)$ such that for all $n \geq 0$ we have
$$
z_i(2^e n + r) = A_i(e,r) z_i(n) + B_i(e,r) z_i(n+1).
$$
Taking $n = 0$ yields $B_i(e,r) = z_i(r)$ (for $i = 1, 2, 3$). Taking $n = 2$, and 
using that $z_1(2) = 1$, $z_3(2) = -1$, and $z_1(3) = z_3(3) = 0$, we get
$A_1(e,r) = z_1(2^{e+1} + r)$ and $A_3(e,r) = -z_3(2^{e+1} + r)$.
Now taking $n=1$ yields $A_2(e,r) = z_2(2^e + r) - z_2(r)$. An immediate induction 
on $e$ proves that, for $r \in [0, 2^e]$, one has
$z_2(2^e + r) - z_2(r) = -z_2(5.2^e + r)$. Hence we can state the following theorem.

\begin{theorem}
Let $(z_1(n))_{n \geq 0}$, $(z_2(n))_{n \geq 0}$, $(z_3(n))_{n \geq 0}$ be the
sequences defined above. Let $e \geq 0$ and $r \in [0, 2^e]$. Then, for all 
$n \geq 0$ we have
$$
\begin{array}{lll}
z_1(2^e n + r) &=& z_1(2^{e+1} + r) z_1(n) + z_1(r) z_1(n+1) \\
z_2(2^e n + r) &=& -z_2(5.2^e n + r) z_2(n) + z_2(r) z_2(n+1) \\
z_3(2^e n + r) &=& -z_3(2^{e+1} + r) z_3(n) + z_3(r) z_3(n+1) \\
\end{array}
$$
and
$$
z_1(2^{e+1}+r)z_1(2n+5) + z_1(r)z_1(2n+3) = -z_1(2^e(n+2)+r) + z_1(2^e(n+1)+r)
$$
$$
-z_2(5.2^e+r)z_2(2n+5) + z_2(r)z_2(2n+3) = -z_2(2^e(n+2)+r) + z_2(2^e(n+1)+r) 
$$
$$
-z_3(2^{e+1}+r)z_3(2n+5) + z_3(r)z_3(2n+3) = z_3(2^e(n+2)+r) + z_3(2^e(n+1)+r). 
$$
\end{theorem}

\subsection{Block-complexity of the Thue-Morse sequence}

Other sequences satisfy the hypotheses of Theorem~\ref{general}, e.g., sequence 
A145865	in \cite{oeis}. An example that we would like to mention is the sequence
$(y(n))_{n \geq 0} = (A005942(n+1))_{n \geq 0}$ with the notation of \cite{oeis}. 
The sequence $(A005942(n))_{n \geq 0}$ is the (block-)complexity of the Thue-Morse
sequence (the Thue-Morse sequence is the fixed point beginning with $0$ of the
morphism $0 \to 01$, $1 \to 10$, see, e.g., \cite{ubi}; its block-complexity is 
the number of distinct factors (blocks) of each length occurring in that sequence). 
It satisfies $A005942(2n)=A005942(n)+A005942(n+1)$, and $A005942(2n+1)=2A005942(n+1)$ 
if $n \geq 2$ (see \cite{Brlek, LucVar}). Hence the sequence $(y(n))_{n \geq 0}$ 
satisfies the hypotheses of Theorem~\ref{general} with $n_0 = 2$, $a=2$, $b=c=1$. 
Note that $y(0) = 2$ and $y(1) = 4$.

\begin{remark}
{\rm The sequence $(A006165(n))_{n \geq 0}$ satisfies the same recurrence properties
as the sequence $(y(n))_{n \geq 0}$ above, but is equal to $1$ for $n=1$ and
$n=2$. As indicated in \cite{oeis} this sequence is related to the Josephus problem.
}
\end{remark}

\section{Final remarks}

For sequences $(z(n))_{n \geq 0}$ satisfying the hypotheses of Theorem~\ref{general},
any subsequence of the form $(z(2^e n + r))_{n \geq 0}$ with $e \geq 0$ and
$r \in [0, 2^e]$ is a linear expression in $(z(n))_{n \geq 0}$ and $(z(n+1))_{n \geq 0}$
for $n \geq n_0$ with coefficients depending on $r$ and $e$ only: this proves the 
$2$-regularity of these sequences (see \cite{allsha}).

Also note that, as visible in the proof of Theorem~\ref{general} above, several 
other relations can be found between the terms of sequences satisfying the hypotheses 
of that theorem.

\section{Acknowledgments} We thank warmly R. Bacher for his comments on a previous version
of this paper.


\begin{thebibliography}{80}

\bibitem{allsha} J.-P. Allouche, J. Shallit, The ring of $k$-regular sequences, 
{\it Theoret. Comput. Sci.} {\bf 98} (1992) 163--197. 

\bibitem{ubi} J.-P. Allouche, J. Shallit, The ubiquitous Prouhet-Thue-Morse sequence,
in Sequences and their applications, Proceedings of SETA'98, C. Ding, T. Helleseth 
and H. Niederreiter (Eds.), 1999, Springer Verlag, 1--16.

\bibitem{Bacher} R. Bacher, Twisting the Stern sequence, Preprint (2010), available
electronically at {\tt http://arxiv.org/abs/1005.5627}

\bibitem{Brlek} S. Brlek, Enumeration of factors in the Thue-Morse word, 
{\it Discrete Applied Math.} {\bf 24} (1989) 83--96.

\bibitem{Carlitz} L. Carlitz, Single variable Bell polynomials, {\it Collect. Math.}
{\bf 14} (1962) 13--25.

\bibitem{Coons-conj} M. Coons, On some conjectures concerning Stern's sequence 
and its twist, {\it Integers} {\bf 11} (2011) \#A35, available
electronically at {\tt http://www.integers-ejcnt.org/vol11.html}

\bibitem{Coons} M. Coons, A correlation identity for Stern's sequence,
{\it Integers\,} {\bf 12} (2012) \#A3, available
electronically at {\tt http://www.integers-ejcnt.org/vol12.html}

\bibitem{DS} K. Dilcher, K. B. Stolarsky, A polynomial analogue to the Stern sequence,
{\it Int. J. Number Theory} {\bf 3} (2007) 85--103. 

\bibitem{LucVar} A. De Luca, S. Varricchio, Some combinatorial properties of the 
Thue-Morse sequence and a problem in semigroups, {\it Theoret. Comput. Sci.}
{\bf 63} (1989) 333--348. 

\bibitem{Northshield} S. Northshield, Stern's diatomic sequence 
0, 1, 1, 2, 1, 3, 2, 3, 1, 4,~..., {\it Amer. Math. Monthly\,} 
{\bf 117} (2010) 581--598.

\bibitem{oeis} {\it The On-Line Encyclopedia of Integer Sequences} available electronically at
{\tt https://oeis.org/}

\bibitem{Reznick} B. Reznick, Regularity properties of the Stern enumeration of the rationals,
{\it J. Integer Seq.} {\bf 11} (2008), Article 08.4.1, available electronically at
{\tt http://www.cs.uwaterloo.ca/journals/JIS/VOL11/Reznick/reznick4.html}

\bibitem{Urbiha} I. Urbiha, Some properties of a function studied by De Rham, 
Carlitz and Dijkstra and its relation to the (Eisenstein-)Stern's diatomic sequence, 
{\it Math. Commun.} {\bf 6} (2001) 181--198.

\end{thebibliography}
\end{document}